\documentclass[11pt,a4paper]{article}
\usepackage{graphics}
\usepackage{color}
\usepackage{graphicx}
\usepackage{graphics}
\usepackage{makeidx}
\usepackage{showidx}
\usepackage{latexsym}
\usepackage{amssymb}
\usepackage{verbatim}
\usepackage{amsmath}
\usepackage{amsthm}
\usepackage{amsfonts}
\usepackage{amssymb,amsmath}
\usepackage{latexsym,amsthm,amscd}
\usepackage[all]{xy}
\newtheorem{thm}{Theorem}
\newtheorem{prop}[thm]{Proposition}
\newtheorem{defn}{Definition}
\newcounter{alphthm}
\newtheorem{lem}[thm]{Lemma}

\newcommand{\be}{\begin{equation}}
\newcommand{\ee}{\end{equation}}
\newcommand{\ben}{\begin{enumerate}}
\newcommand{\een}{\end{enumerate}}
\newcommand{\beq}{\begin{eqnarray}}
\newcommand{\eeq}{\end{eqnarray}}
\newcommand{\beqn}{\begin{eqnarray*}}
\newcommand{\eeqn}{\end{eqnarray*}}

\newcommand{\bpf}{\begin{proof}}
\newcommand{\epf}{\end{proof}}
\newcommand{\bl}{\begin{lem}}
\newcommand{\el}{\end{lem}}
\newcommand{\bp}{\begin{prop}}
\newcommand{\ep}{\end{prop}}
\newcommand{\bd}{\begin{defn}}
\newcommand{\ed}{\end{defn}}
\newcommand{\bt}{\begin{thm}}
\newcommand{\et}{\end{thm}}

\newcommand\bpr{\begin{prop}}
\newcommand\epr{\end{prop}}

\begin{document}
\title{On  the characteristic of projectively invariant pseudo-distance on  Finsler spaces}
 \author{M. Sepasi and  B. Bidabad \footnote{Corresponding author}\\
\small{ Faculty of  Mathematics and  Computer Science,}\\\small{ Amirkabir University  of Technology (Tehran Polytechnic),}\\\small{ Hafez  Ave.,  15914  Tehran, Iran.}\\ \small{ m\_sepasi@aut.ac.ir; bidabad@aut.ac.ir}}
 \date{}
\maketitle

\begin{abstract}
A projective parameter of a geodesic on a Finsler space  is defined to be solution of a certain ODE. Using projective parameter and Funk metric, one can construct a projectively invariant intrinsic pseudo-distance on a Finsler space. In the present work, solutions of the projective parameter's ODE are characterized with respect to the sign of parallel Ricci tensor of a Finsler space. It is shown that the pseudo-distance is trivial on complete Finsler spaces of positive semi-definite Ricci tensor    and it is  a distance on Finsler spaces of parallel negative definite Ricci tensor. These results generalize some results of Kobayashi and Sasaki to Finsler geometry.
\end{abstract}
\section{Introduction}
Let $(M,F)$ and $(M,\bar{F})$ be two Finsler spaces.   If any geodesic of $(M,F)$ coincides with a geodesic of $(M,\bar{F})$ as a set of points and vice versa, then  $F$ and $\bar{F}$ are said to be projectively related.
 It is well-known that a Finsler space $(M,F)$  is projective to another Finsler space $(M,\bar{F})$, if and only if there exists a  1-homogeneous scalar field $P(x,y)$ satisfying $\bar{G}^i(x , y)=G^i(x , y)+P(x,y)y^i$, where $G^i$ and ${\bar{G}}^i$ are corresponding spray vector fields. The scalar field $P(x , y)$ is called the \emph{projective factor}. Let $\gamma:=x^i(t)$ be a geodesic on $(M,F)$. In general, the parameter ``$t$" of $\gamma$, does not remain invariant under projective changes. There is a  parameter  which remains invariant under projective changes called  \emph{projective }parameter.
In Refs. \cite{B, E, T} the projective parameter is defined  for geodesics of general affine connections. In  Ref. \cite{SB1}
 it is carefully spelled out for geodesics of Finsler metrics  as a solution of the following ODE
\begin{equation}\label{e16}
\{p,s\}:= \frac{\frac{d^3p}{ds^3}}{\frac{dp}{ds}}-\frac{3}{2}\big[\frac{\frac{d^2p}{ds^2}}{\frac{dp}{ds}}\big]^2 =\frac{2}{n-1}F^2Ric = \frac{2}{n-1} Ric_{jk}\frac{d{x}^j}{ds}\frac{d{x}^k}{ds},
\end{equation}
where $\{p,s\}$ is  known  in the literature  as \emph{Schwarzian derivative}
 and  ``$s$" is the arc length parameter of $\gamma$.  The projective parameter is unique up to all linear fractional transformations
 \begin{equation}\label{PropertyOfSchwarzian}
 \{ \frac{a p + b}{c p +d} , s \} = \{p , s\},
 \end{equation}
where, $ad-bc\neq 0$.
 Previously, the present authors, using projective parameter and Funk metric on the open interval $(-1,1)$, studied  an intrinsic  projectively invariant pseudo-distance  denoted by $d_M$, cf., \cite{SB1}. Next, in \cite{SB2} it is shown that in a complete Einstein Finsler space with negative constant Ricci scalar, the intrinsic projectively invariant pseudo-distance is a constant multiple of the Finslerian distance. Therefore, as a corollary, it is deduced that two  projectively  related complete Finsler spaces with constant negative Ricci scalar  are homothetic. The last result is previously obtained by Z. Shen using another technic of proof. See Ref. \cite{S}.\\In the present work, we investigate the differential equation (\ref{e16}) when
the Ricci tensor is parallel with respect to any of Berwald, Chern or Cartan
connection and will present the solution. More precisely, we prove
\bt
Let $(M,F)$ be a Finsler space of parallel Ricci tensor. Then the Ricci tensor is constant along geodesics parameterized by arc-length, and solutions of (\ref{e16}) are given as follows.
\begin{itemize}
\item[i)] If $\{p,s\}=c^2$ with $c>0$ then
\begin{equation}\label{p1}
p=\frac{\alpha cos(cs)+\beta sin(cs)}{\gamma cos(cs)+\delta sin(cs)}.
\end{equation}
\item[ii)] If $\{p,s\}=-c^2$ with $c>0$ then
\begin{equation}\label{p2}
p=\frac{\alpha e^{cs}+\beta e^{-cs}}{\gamma e^{cs}+\delta e^{-cs}}.
\end{equation}
\item[iii)] If $\{p,s\}=0$ then
\begin{equation}\label{p3}
p=\frac{\alpha+\beta s}{\gamma+\delta s}.
 \end{equation}
 \end{itemize}
 \et
Here, a new approach to the study of the intrinsic pseudo-distance is considered and following results are obtained.
\bt
Let $(M,F)$ be a connected complete Finsler space of positive semi-definite  Ricci tensor. Then the intrinsic projectively invariant pseudo-distance is trivial, that is $d_M=0$.
\et
\bt
Let $(M,F)$ be a connected (complete) Finsler space of negative-definite parallel Ricci tensor in  Berwald or Chern connection. Then  the intrinsic projectively invariant pseudo-distance, $d_M$, is a (complete) distance.
\et
 These Theorems are generalizations of some results in \cite{K} and \cite{KS}.
\section{Preliminaries }

Here and every where in this work the differential  manifold $M$ is supposed to be a connected differential manifold.
A (globally defined) Finsler structure on a differential manifold $M$ is a function $F: TM\rightarrow [0 , \infty) $ with the  properties, i) Regularity: $F$ is $C^{\infty}$ on the entire slit tangent bundle $TM_0$,
ii) Positive homogeneity:  $F(x , \lambda y) = \lambda F(x , y)$ for all $\lambda > 0$,
iii) Strong convexity: The Hessian matrix $(g_{ij}) := ({[1/2F^2]}_{y^iy^j})$ is positive-definite at every point of $TM_0$. The pair $(M,F)$ is known as a Finsler space.

Every Finsler structure $F$ induces a spray $\textbf{G}=y^i\frac{\partial}{\partial x^i}-G^i(x,y)\frac{\partial}{\partial y^i}$ on $TM$, where $G^i(x,y):=\frac{1}{2} g^{il}\{[F^2]_{x^k y^l}y^k-[F^2]_{x^l}\}$. $\textbf{G}$ is a globally defined vector field on $TM$.  Projection of a flow line of $\textbf{G}$ on $M$
is called a geodesic . Differential equation of a geodesic in local coordinate is given by $\frac{d^2x^i}{ds^2}+G^i(x(s),\frac{dx}{ds})=0$, where $s(t) = \int_{t_0}^{t} F(\gamma , \frac{d\gamma}{dr}) dr$ is the arc length parameter.

For a  non null $y \in T_xM$, the Riemann curvature
$\textbf{R}_y : T_xM \rightarrow T_xM$ is defined by $\textbf{R}_y(u)=R^i_k u^k \frac{\partial}{\partial x^i}$, where $R^i_k(y):=\frac{\partial G^i}{\partial x^k}-1/2\frac{{\partial}^2G^i}{\partial y^k \partial x^j}y^j+G^j \frac{{\partial}^2G^i}{\partial y^k \partial y^j}-1/2\frac{\partial G^i}{\partial y^j}\frac{\partial G^j}{\partial y^k}$.
The \emph{Ricci Scalar} is defined by $Ric:= {R^i}_i$. \cite{BCS}.
 In the present work, we use the definition of  \emph{Ricci tensor}
introduced  by Akbar-Zadeh,  as follows  $Ric_{ik} :=\frac{1}{2}(F^2Ric)_{y^iy^k}$. cf., \cite{A}.
Moreover, by homogeneity we have $Ric_{ik} {\ell}^i{\ell}^k = Ric$.

Let $G^{i}:={\gamma}^i_{jk}y^iy^j$, where ${\gamma^i}_{jk}:=1/2g^{is}(\frac{\partial g_{sj}}{\partial x^k}-\frac{\partial g_{jk}}{\partial x^s}+\frac{\partial g_{ks}}{\partial x^j})$, ${N^i}_j :=1/2 \frac{\partial G^i}{\partial y^j}$, $l^i:=\frac{y^i}{F}$, and  ${\widehat{l}}:=l^i\frac{\delta}{\delta x^i}=l^i( \frac{\partial}{\partial x^i}-{N^k}_i \frac{\partial}{\partial y^k})$.  See Ref. \cite{BCS}.

\section{Projective parameter for Ricci parallel Finsler spaces}

Let the Ricci tensor of $(M,F)$  be parallel with respect to any of  Cartan, Berwald or Chern connection. We recall
 the \emph{Abel's identity} in ordinary differential equations as follows.

\emph{Consider the second -order linear ordinary differential equation;}
\begin{equation}\label{abel}
y^{''}+P(x)y^{'}+Q(x)y=0.
\end{equation}
\emph{call the two linearly independent solutions, $y_1(x)$ and $y_2(x)$.Then, the Wronskian of $y_1$ and $y_2$,} $w(y_1,y_2)=y_1y^{'}_2- y_2y^{'}_1 $ \emph{satisfies} $w^{'}+pw=0$, \emph{therefore}
\begin{equation}\label{abel1}
w=w_0 e^{-\int P(x)dx}.
\ee
\setcounter{thm}{0}
\bp
If  $y_1$ and $y_2$ are linearly independent solution of  the ordinary differential equation
\begin{equation}\label{a}
y^{''}+Q(s)y(s)=0,
\ee
where $Q(s)=\frac{1}{n-1}Ric_{ij}\frac{dx^i}{ds}\frac{dx^j}{ds}$, then the general solution of (\ref{e16}) is given by
\begin{equation}\label{b}
u(t)=\frac{\alpha y_1 + \beta y_2}{\gamma y_1 + \delta y_2},
\ee
with $\alpha \delta-\beta \gamma\neq 0$
\ep
\begin{proof}
According to (\ref{PropertyOfSchwarzian}), it saffices to show that ${y_1}/{y_2}$ is a solution of (\ref{e16}). $P(x)$ in (\ref{abel1}) is zero, so the Wronskian $w(y_1,y_2)$ is constant. We
may assume that $w(y_1,y_2)=1$. Then, $u^{'}={1}/{y^2_2}$. So that ${u^{''}}/{u^{'}}= -2{y^{'}_2}/{y_2}$, and
$$(\frac{u^{''}}{u^{'}})^{'}=\frac{-2y^{''}_2 y_2 +2(y^{'}_2)^2}{y^2_2}=-2\frac{y^{''}_2 }{y^{''}_2 }+2(\frac{y^{'}_2}{y_2})^2,$$
$$\frac{u^{'''}u^{'}-(u^{''})^2}{(u^{'})^2}=-2\frac{(-Q(s))y_2(s)}{y_2(s)}+\frac{1}{2}(\frac{u^{''}}{u^{'}})^{2},$$
$$\frac{u^{'''}}{u^{'}}-(\frac{u^{''}}{u^{'}})^{2}=2Q(s)+\frac{1}{2}(\frac{u^{''}}{u^{'}})^{2},$$
$$\frac{u^{'''}}{u^{'}}-\frac{3}{2}(\frac{u^{''}}{u^{'}})^{2}=2Q(s).$$
This completes the proof.
\end{proof}
\emph{Proof of Thoerem 1.}
  Let  the Ricci tensor be parallel with respect to Cartan connection. We denote the horizontal and vertical Cartan covariant derivative of Ricci tensor by
   $\bigtriangledown^c_{\frac{\delta}{\delta x^k}}Ric_{ij}$ and $\bigtriangledown^c_\frac{\partial}{\partial y^k}Ric_{ij}$ respectively. we have
 \begin{equation}\label{c}
\bigtriangledown^c_{\frac{\delta}{\delta x^k}}Ric_{ij}=\frac{\delta Ric_{ij}}{\delta x^k}- Ric_{ir}{\Gamma^r}_{jk}-Ric_{jr}{\Gamma^r}_{ik}=0,
\ee
\begin{equation}\label{d}
\bigtriangledown^c_{\frac{\partial}{\partial y^k}}Ric_{ij}=\frac{\partial Ric_{ij}}{\partial y^k}- Ric_{ir}\frac{{A^r}_{jk}}{F}-Ric_{jr}\frac{{A^r}_{ik}}{F}=0,
\ee
where ${\Gamma^i}_{jk}=\frac{1}{2}g^{ih}(\frac{\delta g_{hj}}{\delta x^k}+\frac{\delta g_{kh}}{\delta x^j}-\frac{\delta g_{jk}}{\delta x^h})$  and ${A^i}_{jk}:=g^{ih}A_{hjk}=g^{ih}\frac{F}{4}\frac{\partial{g_{ij}}}{\partial y^k}$ is the coefficient of  \emph{Cartan} tensor.
Consider  the geodesic $\gamma:=x^i(s)$, where ``s" is the arc-length parameter. Contracting (\ref{c}) by $\frac{dx^i}{ds}\frac{dx^j}{ds}\frac{dx^k}{ds}$ gives $$\frac{dx^i}{ds}\frac{dx^j}{ds}\frac{dx^k}{ds}(\frac{\partial Ric_{ij}}{\partial x^k}-{N^l}_k\frac{\partial Ric_{ij}}{\partial y^l})-$$
$$\frac{dx^i}{ds}\frac{dx^j}{ds}\frac{dx^k}{ds}(Ric_{ir}{\Gamma^r}_{jk})-\frac{dx^i}{ds}\frac{dx^j}{ds}\frac{dx^k}{ds}(Ric_{jr}{\Gamma^r}_{ik})=0.$$
Using (\ref{d}), and the property $y^j {A^i}_{jk}=y^k {A^i}_{jk}=0$ of Cartan tensor, we have
$$\frac{dx^i}{ds}\frac{dx^j}{ds}\frac{dRic_{ij}}{ds}-\frac{dx^i}{ds}\frac{dx^j}{ds}\frac{dx^k}{ds}{N^l}_k(Ric_{ir}\frac{{A^r}_{jl}}{F}+Ric_{jr}\frac{{A^r}_{il}}{F})$$
$$-2\frac{dx^i}{ds}\frac{dx^j}{ds}\frac{dx^k}{ds}Ric_{jr}{\Gamma^r}_{ik}=0.$$
Therefore
$$\frac{dRic_{ij}\frac{dx^i}{ds}\frac{dx^j}{ds}}{ds}-2Ric_{ij}\frac{d^2x^i}{ds}\frac{dx^j}{ds}-0+2Ric_{rj}\frac{d^2x^r}{ds}\frac{dx^j}{ds}=0,$$
and
\begin{equation}\label{e}
Ric_{ij}\frac{dx^i}{ds}\frac{dx^j}{ds}=constant.
\ee
Following the  method  just used, we can prove that if the Ricci tensor is parallel with respect to the Berwald or Chern connection then along  the geodesic $\gamma$ parameterized by arc-length,  we have $Ric_{ij}\frac{dx^i}{ds}\frac{dx^j}{ds}=constant$.
\\Considering the above assertion and Lemma 1, the equation (\ref{e16})  reduces to a second order ODE with constant coefficient. Thus with respect to the sign of Ricci tensor, one can explicitly determine   a projective parameter ``$p$" as an elementary function of ``s" by (\ref{p1}), (\ref{p2}) and (\ref{p3}). This completes the proof.
$$\ \ \ \ \ \ \ \ \ \ \ \ \ \ \ \ \ \ \ \ \ \ \ \ \ \ \ \ \ \ \ \ \ \ \ \ \ \ \ \ \ \ \ \ \ \ \ \ \ \ \ \ \ \ \ \ \ \ \ \ \ \ \ \ \ \ \ \ \ \ \ \ \ \ \ \ \ \ \ \ \ \ \ \ \ \ \ \ \ \ \ \ \ \ \ \ \ \ \ \ \ \ \ \ \ \ \ \ \ \ \ \ \ \ \ \Box$$
\section{Positive semi-definite Ricci tensor }
 Let consider  the Funk metric $L_f$ and the Funk distance $D_f$ on $I$ by
 \begin{equation}\label{e20}
L_f=\frac{1}{k} (\frac{\mid y \mid}{1-u^2}+\frac{uy}{1-u^2}),
 \end{equation}
\begin{equation}\label{e21}
D_f(a,b) = \frac{1}{2k} (\mid \ln \frac{(1-a)(1+b)}{(1-b)(1+a)} \mid + \ln \frac{(1-a^2)}{(1-b^2)}) \quad a,b\in I.
\ee
See Refs. \cite{O,SB1} for a survey.
Let $f(u)$ be a geodesic on $(M,F)$. If $u$ is a projective parameter then $f$ is said to be \emph{projective}.

Given any two points $x$ and $y$ in $(M,F)$, we consider a chain $\alpha$ of geodesic segments joining these points. That is;\\
 i)a chain of points $x = x_0 , x_1 , ... ,x_k = y$ on $M$; \\ii)pairs of points $a_1,b_1 ,..., a_k,b_k$ in $I$; \\iii)projective maps $f_1,...,f_k$, $f_i: I \rightarrow M $ such that
$f_i(a_i) = x_{i-1}, f_i(b_i) = x_i, \quad i = 1,...,k$.\\
By virtue of the Funk distance $D_f(.,.)$ on $I$ we define the length $L(\alpha)$ of the chain $\alpha$ by
$L(\alpha):= \Sigma_i D_f(a_i , b_i)$, and we put
\begin{equation}\label{e19}
d_M(x , y):= inf L(\alpha),
\end{equation}
where the infimum is taken over all  the chains $\alpha$ of geodesic segments from $x$ to $y$.
One can easily prove the following Lemma.
\setcounter{thm}{0}
\bl
Let $(M, F)$ be a Finsler space. Then for any points $x$, $y$, and $z$ in $M$,  $d_M$ satisfies
 \begin{itemize}
\item[ i)] $d_M(x , y) \neq d_M(y , x)$;
 \item [ii)]$d_M(x , z) \leq d_M(x , y) + d_M(y , z)$;
 \item [iii)]If $x = y$ then $d_M(x , y) = 0$ but the inverse is not  always true.
\end{itemize}
\el
Traditionally, call $d_M(x,y)$ the \emph{ pseudo-distance} of any two points $x$ and $y$ on $M$.
From the property (\ref{PropertyOfSchwarzian}) of Schwarzian derivative, and the fact that the projective parameter is invariant under fractional transformation, the pseudo-distance $d_M$ is projectively invariant.
\subsection{Proof of Theorem 2}
In this section we bring first some Lemmas which will be used in the proof of Theorem 2.
\bl
Let $(M,F)$ be a complete Finsler space. Consider $x_0$ and $x_1$ on $M$. If there exists a geodesic $x(u)$ with projective parameter $u$, $-\infty <u <+ \infty$, such that $x_0=x(u_0)$ and $x_1=x(u_1)$ for some $u_0$ and $u_1$ in $\mathbb{R}$ then
$$d_M(x_0,x_1)=0$$
\el
\begin{proof}
Let us denote the linear equation of the segment passing through the points $(u_0,-1/2)$ and $(u_1,1/2)$  $\widehat{u}={u}/({u_1-u_0})-{1}/{2}({u_1+u_0)}/(u_1-u_0)$. $\widehat{u}$ is a linear transformation of $u$ and is also a projective parameter. We have $-\frac{1}{2}<\widehat{u}<\frac{1}{2}$ when $u_0<u<u_1$.
Next, we consider  the chain $\alpha$ of projective maps,  $a_n$ and $b_n$ where
$$f_n=x(n\widehat{u}) \qquad a_n=-\frac{1}{2n}, \qquad b_n=\frac{1}{2n}.$$
We note that $f_n(-\frac{1}{2n})=x(n(-\frac{1}{2}n))=x(-\frac{1}{2})=x(u_0)$ and $D_f(-\frac{1}{2n},\frac{1}{2n}) = \frac{1}{2k} (\mid \ln \frac{(1+\frac{1}{2n})(1+\frac{1}{2n})}{(1-\frac{1}{2n})(1-\frac{1}{2n})} \mid + \ln \frac{(1-\frac{1}{4n^2})}{(1-\frac{1}{4n^2})})$. Considering $n$ sufficiently large, we have $d_M(x_0,x_1)=infL(\alpha)=0$. This completes the proof.
\end{proof}
\bl
Let $(M,F)$ be a complete Finsler space and $x(s)$ be a geodesic with arc-length parameter $-\infty<s<\infty$. Assume that there exists a (finite or infinite) sequence of open intervals $I_i=(a_i,b_i)$, $i=0,\pm 1,\pm 2,...$, such that;\\ i) $a_{i+1}\leq b_i$, $lim_{i\rightarrow -\infty} a_i=-\infty$ and $lim_{i\rightarrow \infty}b_i=+\infty$  So that $\bigcup_i \overline{I_i}=(-\infty,+\infty)$;\\ ii) in each interval $I_i=(a_i,b_i)$, a projective parameter ``$u$" moves from $-\infty$ to $+\infty$ whenever $t$ moves from $a_i$ to $b_i$. Then, for any pair of points $x_0$ and $x_1$ on this geodesic, we have
$$d_M(x_0,x_1)=0.$$
\el
\begin{proof}
By Lemma 1, the distance between any two points in the same interval $I_i$ is zero.  Two consecutive open intervals $I_i$ and $I_{i+1}$ have either a point  as a boundary point or an interval in common. In either case, given $\epsilon >0$, there exist  the points $s_i$ and $s_{i+1}$ in $I_i$ and $I_{i+1}$ respectively such that $d(x(s_i),x(s_{i+1})) <\epsilon$. This completes the proof.
\end{proof}
The following Lemmas help  to construct open intervals $I_i$ as in Lemma 3. See the proofs in Ref. \cite{KS}
\bl
In the ODE (\ref{a}), if $Q(s)\geq 0$ for all $s \in \mathbb{R}$ then every solution $y(s)$  has at least one zero unless $Q(s)=0$ and $y(s)$ is constant $c\neq 0$.
\el
It is worth noting to recall the \emph{Sturm's separation} theorem as follows:\\
\emph{Given a homogeneous second order linear differential equation and two continuous linear independent solutions $v(x)$ and $u(x)$ with $x_0$ and $x_1$ successive  roots of $v(x)$ then $u(x)$ has exactly one root in the open interval $(x_0,x_1)$}.
\bl
Let $y_1(s)$ and $y_2(s)$ be two linearly independent solutions of (\ref{a}). If $a$ and $b$ are two consecutive zeros of $y_2(s)$ then $u={y_1(s)}/{y_2(s)}$ or $u=-{y_1(s)}/{y_2(s)}$ is a projective parameter in interval $(a,b)$ which moves from $-\infty$ to $+\infty$ as $s$ moves from $a$ to $b$.
\el

The differential equation (\ref{a}) is said to be be oscillatory at $s=\pm \infty$  if the zeros
$$...<a_{-2}<a_{-1}<a_{0}<a_1<a_2<...$$
of the solution $y(s)$ have the property that $lim_{h\rightarrow -\infty}a_h=-\infty$ and $lim_{k\rightarrow+\infty}=+\infty$. Then the sequence of intervals $I_i=(a_i,a_{i+1})$ satisfies the condition of Lemma 3. This fact  proves Theorem 2 in this case. \\We consider the case   the differential equation (\ref{a}) is nonoscillatory at $s=+\infty$. That is, $y_2(s)$ does not vanish for sufficiently large $s$. According to Sturm's theorem, this condition is independent of  choice of a particular solution $y_2(s)$.
\bl
If the differential equation (\ref{a}) is nonoscillatory at $s=+\infty$, then there is a solution $y_2(s)$ which is uniquely determind up to a constant factor satisfying
 \begin{equation}\label{e}
lim_{s\rightarrow +\infty}\frac{y_2(s)}{y_1(s)}=0,
\end{equation}
for any solution $y_1(s)$ linearly independent of $y_2(s)$.
\el
A solution $y_2(s)$ in Lemma 6 is called a \emph{principal} solution.
Here, we consider a weaker version of \emph{Comparison Theorem of Sturm} as follows
\bl
Consider two differential equations
$$(i) y^{''}(s)+Q_1(s)y(s)=0,\qquad (ii)  y^{''}(s)+Q_2(s)y(s)=0,$$
with  $Q_1(s)\geq Q_2(s)$. Let $y_1(s)$ and $y_2(s)$ be solutions of (i) and (ii) respectively such that
\begin{equation}\label{g}
\frac{y_1^{'}(a)}{y_1(a)}\leq \frac{y_2^{'}(a)}{y_2(a)}.
\end{equation}
 If $y_1(s)$ and $y_2(s)$ have no zero in the interval $a<s<+\infty$, then for $s>a$
\begin{equation}\label{h}
\frac{y_1^{'}(s)}{y_1(s)}\leq \frac{y_2^{'}(s)}{y_2(s)}.
\end{equation}
If $y_2(a)=0$,  then the term $y_2^{'}(a)/y_2(a)$ is considered to  be $\infty$.
\el
See Refs.\cite{DK} and \cite{KS} for more details.
\bl
Assume that the differential equation (\ref{a}) is nonoscillatory at $s=+\infty$ and that $Q(s)\geq 0$. Let
$y_(s)$ be a principal solution as in Lemma 6. If  $a$ is the largest zero of $y_2(s)$ and if $y_1(s)$ is a solution linearly independent of $y_2(s)$, then $y_1(s)$ vanishes at some $s>a$.
\el
We are now in a position to complete the proof of the theorem 2 where the differential equation (\ref{a}) is nonoscillatory at $s=+\infty$ or $s=-\infty$.
\\If (\ref{a}) is nonoscillatory at $s=+\infty$ but oscillatory at $s=-\infty$, we take a principal solution $y_2(s)$ and another solution  $y_1(s)$ linearly independent of $y_2(s)$. Let $...<a_{-2}<a_{-1}<a_0<a_1<a_2<...<a_k$ be the zeros of $y_2(s)$. Then the sequence of intervals $I_i=(a_i,a_{i+1})$, for $i=...,-2,-1,0,1,2,...,k$ with $a_{k+1}=+\infty$, equipped with a projective parameter $u=y_1/y_2$ or $u=-y_1/y_2$ satisfy the requirements of Lemma 3. We note that Lemma 8 implies  that $u$ is a projective parameter in the last interval $I_k=(a_k,+\infty)$. If (\ref{a}) is nonoscillatory at $s=-\infty$ but oscillatory at  $s=+\infty$, we replace Lemma 6 and Lemma 8 by the analogous Lemmas for $s=-\infty$. Assume that (\ref{a}) is nonoscillatory at $\pm \infty$. Let $y_2(s)$ be  a principal solution for $s=+\infty$ and not for $s=-\infty$. Let $y_1(s)$ be a principal solution for $s=-\infty$ then $y_1(s)$ and $y_2(s)$ are linearly independent. We obtain a sequence of intervals $I_i$, $i=0,1,...,k$ with a projective parameter $u=y_1/y_2$, $-y_1/y_2$, $y_2/y_1$ or $-y_2/y_1$ satisfying  the requirements of Lemma 3. In this case, there are some overlaps among these intervals. \\If $y_2(s)$ is a principal solution for both $s=+\infty$ and $s=-\infty$ then we consider $y_1(s)$ as a solution linearly independent of $y_2(s)$. We obtain a sequence of intervals  $I_i$, $i=0,1,...,k$, with a projective parameter $u=y_1/y_2$ or $-y_1/y_2$ satisfying the requirements of Lemma 3. In this case, there are no overlaps of intervals. This completes the proof of Theorem 2.
$$\ \ \ \ \ \ \ \ \ \ \ \ \ \ \ \ \ \ \ \ \ \ \ \ \ \ \ \ \ \ \ \ \ \ \ \ \ \ \ \ \ \ \ \ \ \ \ \ \ \ \ \ \ \ \ \ \ \ \ \ \ \ \ \ \ \ \ \ \ \ \ \ \ \ \ \ \ \ \ \ \ \ \ \ \ \ \ \ \ \ \ \ \ \ \ \ \ \ \ \ \ \ \ \ \ \ \ \ \ \ \ \ \ \ \ \Box$$
\section{Parallel negative-definite Ricci tensor  }
We recall the following theorem which will be used in the sequel. \\
\textbf{Theorem A.} \cite{SB1}
\emph{Let $(M,F)$ be  a connected (complete) Finsler space for which the Ricci tensor satisfies $Ric_{ij}\leq g_{jk}$, as matrices, for a positive constant $c$. Then $d_M$ is a (complete) distance.}\\
\emph{Proof of Theorem 3.}
Consider the Finsler structure  $\hat{F}(x,y)=\sqrt{-Ric_{ij}(x,y)y^iy^j}$. It suffices to show  that the spray coefficients of  $\hat{F}$ and $F$ are equal, that is ${\hat{G}}^i=G^i$. According to   Theorem A,  $d_M$ is a (complete) distance.\\
We have
\begin{eqnarray}\label{G^i}
&&{\hat{G}}^i= 1/2 (-Ric)^{ih}(\frac{{\partial}^2 {\hat{F}}^2}{\partial y^h \partial x^j}y^j-\frac{\partial {\hat{F}}^2 }{\partial x^h}) \nonumber \\
&&=1/2 (-Ric)^{ih}(\frac{{\partial}^2 (-Ric_{lr}y^l y^r)}{\partial y^h \partial x^j}y^j-\frac{\partial (-Ric_{lr}y^l y^r)}{\partial x^h}) \nonumber\\
&&=1/2(-Ric)^{ih}(-2\frac{\partial Ric_{hl}y^l}{\partial x^j}y^j+\frac{\partial Ric_{lr}}{\partial x^h}y^l y^r) \nonumber \\
&&=Ric^{ih}\frac{\partial Ric_{hl}y^l}{\partial x^j}y^j- 1/2 Ric^{ih}\frac{\partial Ric_{lr}}{\partial x^h}y^l y^r.
\end{eqnarray}
Let $\bigtriangledown^b$ denote the Berwald connection and Ricci tensor be parallel with respect to the Berwald connection. Similar arguments as follows hold well for Chern connection. We have
 \begin{equation}\label{Ber1}
\bigtriangledown^b_{\frac{\delta}{\delta x^j}}Ric_{hl}=\frac{\delta Ric_{hl}}{\delta x^j}- Ric_{hr}{G^r}_{lj}-Ric_{lr}{G^r}_{hj}=0,\quad {G^r}_{lj}=1/2\frac{\partial^2 G^r}{\partial y^j y^l}.
\ee
\begin{equation}\label{Ber2}
\bigtriangledown^b_{\frac{\partial}{\partial y^k}}Ric_{ij}=\frac{\partial Ric_{ij}}{\partial y^k}=0.
\ee
Contracting (\ref{Ber1}) in $Ric^{ih}y^jy^l$ we have
$$Ric^{ih}y^jy^l\frac{\partial Ric_{hl}}{\partial x^j}-Ric^{ih}Ric_{ha}G^a-1/2Ric^{ih} Ric_{la}\frac{\partial G^a}{\partial y^h}y^l=0.$$
\begin{equation}\label{Ber3}
Ric^{ih}y^jy^l\frac{\partial Ric_{hl}}{\partial x^j}-G^i-1/2Ric^{ih} Ric_{la}\frac{\partial G^a}{\partial y^h}y^l=0.
\ee
On the other hand
$$-1/2 Ric^{ih}y^ry^l\frac{\partial Ric_{lr}}{\partial x^h}+ 1/2Ric^{ih}y^ry^lRic_{la}{G^a}_{rh}+1/2 Ric^{ih}y^ry^lRic_{ra}{G^a}_{lh}=0.$$
\begin{equation}\label{Ber4}
-1/2 Ric^{ih}y^ry^l\frac{\partial Ric_{lr}}{\partial x^h}+ 1/2Ric^{ih}y^r Ric_{ra}\frac{\partial G^a}{\partial y^h}=0.
\ee
Considering (\ref{G^i}), (\ref{Ber3}) and (\ref{Ber4}) we have ${\hat{G}}^i=G^i$. Therefor  $Ric_{ij}$ is the Ricci tensor of $\hat{F}$ too. According to Theorem A, $d_{M}$ is a (complete) distance.
$$\ \ \ \ \ \ \ \ \ \ \ \ \ \ \ \ \ \ \ \ \ \ \ \ \ \ \ \ \ \ \ \ \ \ \ \ \ \ \ \ \ \ \ \ \ \ \ \ \ \ \ \ \ \ \ \ \ \ \ \ \ \ \ \ \ \ \ \ \ \ \ \ \ \ \ \ \ \ \ \ \ \ \ \ \ \ \ \ \ \ \ \ \ \ \ \ \ \ \ \ \ \ \ \ \ \ \ \ \ \ \ \ \ \ \ \Box$$
 
  \end{document}